\documentstyle[epsf,sectioneqno]{article}
\newtheorem{theorem}{Theorem}[section]
\newtheorem{lemma}[theorem]{Lemma}
\newtheorem{cor}[theorem]{Corollary}

\newtheorem{condition}[theorem]{Condition}

\title{Analytic continuation and resonance-free regions for Sturm-Liouville potentials with power decay}
\author{    B.M.Brown, M.S.P. Eastham \\ 
 Department of Computer Science,  University of Cardiff, Cardiff, CF24 3XF,
U.K.  } 
\begin{document}
\maketitle
\newcommand{\beq}{\begin{equation}}
\newcommand{\enq}{\end{equation}}
\newcommand{\Si}{ {{\rm Sin}}}
\bibliographystyle{plain}
 
\def\R{ {\bf R}}
\def\N{ {\bf N}}
  
\section{Introduction}
We consider the Sturm-Liouville equation
\begin{equation}
y^{''}(x)+\{ \lambda -q(x) \} y(x)=0\;\;\;\;\;(0 \leq x < \infty )
\label{eq:1.1}
\end{equation}
with a boundary condition
\begin{equation}
y(0) \cos \alpha + y^{'}(0) \sin \alpha =0,
\label{eq:1.2}
\end{equation}
$\lambda$ being the complex spectral parameter. As usual, $\alpha$ is
real and the potential $q$ is real-valued and locally integrable on
$[0,\infty)$. We further assume throught the paper that $q$ decays as $x
\rightarrow \infty$ in the sense that
\begin{equation}
q \in L(0,\infty).
\label{eq:1.3}
\end{equation}
\par
Let us write $\lambda=z^2$, where $0 \leq {\rm arg \; } z < \pi$ when
$0 \leq {\rm arg \; } \lambda < 2 \pi$. Then (\ref{eq:1.3}) implies that there is
a solution   $\psi(x,z)$ of (\ref{eq:1.1}) such that
\begin{equation}
\psi(x,z) \sim \exp(izx),\;\;\;\;\psi^{'}(x,z) \sim iz \exp(izx)
\label{eq:1.4}
\end{equation}
as $x \rightarrow \infty$, and $\psi(x,z)$ is analytic in $z$ for im $z>0$
\cite[Theorem 1.9.1]{MSPE89}.
Then $\psi(x,z)$ is the Weyl $L^2(0,\infty)$ solution of
(\ref{eq:1.1}) when $\lambda$ is non-real and it forms the basis of the
Weyl-Titchmarsh spectral theory of (\ref{eq:1.1}) \cite[Chapter 9]{CL},
\cite{titc62},\cite{weyl10}. A central result of this spectral theory is the
existence of a spectral function $\rho_\alpha(\mu)\;\;\;(-\infty< \mu <
\infty)$ which is piecewise constant in $(-\infty,0)$ and locally absolutely
continuous in $[0,\infty)$ with $\rho^{'}_\alpha (\mu)>0$
\cite[section 5.7]{titc62}, \cite[p. 264]{weyl10}. In particular,
(\ref{eq:1.4}) leads to the Kodaira formula
\begin{equation}
\pi \rho^{'}_\alpha(\mu)= \mu^{1/2} / \mid
\Psi(\mu^{1/2})\mid^2\;\;\;\;\;(\mu>0)
\label{eq:1.5}
\end{equation}
where 
\begin{equation}
\Psi(z)=\psi(0,z)\cos \alpha + \psi^{'}(0,z)\sin \alpha
\label{eq:1.6}
\end{equation}
\cite[p. 940]{K49}.
Since the only possible eigenvalues of the problem
(\ref{eq:1.1})-(\ref{eq:1.3}) lie on the negative real $\lambda$-axis,
$\Psi(z)$ has no zeros for $0\leq {\rm arg \;} z < \pi$
expect possibly when arg $z=\frac{1}{2}\pi$.
\par
In addition to (\ref{eq:1.5}), the Weyl-Titchmarsh function $m_\alpha(\lambda)$
\cite[Chapter 9]{CL},\cite[Chapter 2]{titc62} also involves $\Psi$ in the form
\begin{equation}
m_\alpha(\lambda)=\{\psi^{'}(0,z)\cos \alpha - \psi(0,z) \sin \alpha \}/\Psi(z)
\label{eq:1.7}
\end{equation}
again with $0 \leq {\rm arg \;} z < \pi$. Now $m_\alpha(\lambda)$ is related to
the
Green's function and to the resolvent operator of (\ref{eq:1.1})-(\ref{eq:1.2})
in the Hilbert space $L^2(0,\infty)$, and the question arises whether these
three spectral objects have analytic continuations into the so-called
unphysical sheet $\pi \leq {\rm arg \; } z <2\pi$.
As far as the Green's function and resolvent are concerned, this 
question can be posed, not only for (\ref{eq:1.1})-(\ref{eq:1.2}),
but also for the corresponding Schr\"{o}dinger equation in  two or more
dimensions. However, in the case of (\ref{eq:1.1}) itself, it is a question of
the analytic continuation of $\psi(0,z)$ and $\psi^{'}(0,z)$ in (\ref{eq:1.6})
and (\ref{eq:1.7}).
\par
Analytic continuation into the strip $-\frac{1}{2}a <{\rm im \;} z <0$ was
established in \cite{DMT66} subject to a strengthening of (\ref{eq:1.3})
to 
\begin{equation}
q(x)=O(e^{-ax})\;\;\;\;(x \rightarrow \infty)
\label{eq:1.8}
\end{equation}
for some $a >0$ (see also \cite[section 2.2]{HS87}), and we refer again to
\cite{DMT66} for a description of earlier work in this direction.
Allowing $a$ to be arbitrarily large in (\ref{eq:1.8}) leads to the class of
super-exponentially decaying potentials for which
\begin{equation}
q(x)=O(e^{-xf(x)})\;\;\;\;(x \rightarrow \infty)
\label{eq:1.9}
\end{equation}
 with some $f(x) \rightarrow \infty$, and then we have analytic continuation
into the whole of im $z<0$. \cite{F97}, \cite{MH99}.
Two other specalisations of (\ref{eq:1.8}) where again there is analytic
continuation into the whole of im $z<0$ are
\begin{displaymath}
q(x)=e^{-ax}p(x)
\end{displaymath}
with $p(x) $ periodic and
\begin{displaymath}
q(x)=({\rm const.})x^Ne^{-ax}
\end{displaymath}
where $N(\geq 1)$ is an integer \cite{BEM99}.
\par
Once analytic continuation has been effected, the possibility is opened up of
$\Psi(z)$
having zeros in the unphysical sheet.  
Such zeros are called resonances and, by (\ref{eq:1.7}), they are singular
spectral points associated with the Green's function and resolvent operator.
For potentials of the class (\ref{eq:1.9}), the asymptotic distribution of
resonances was obtained in \cite{F97} and, by another method, also in
\cite{MH99} along with other results on the location of resonances.
In particular \cite[Theorem 3.8]{MH99}, there is a resonance-free strip $-b
\leq {\rm im \; } z <0$
subject to $q$ having a suitably small norm.
\par
All these existing results require exponential decay of the potential $q$.
In this paper, we allow $q$ to have only power decay $O(x^{-\gamma})$ $(x
\rightarrow \infty)$ for some $\gamma >1$ and,
under further conditions on the analyticity of $q$, we establish analytic
continuation of $\Psi(z)$ into a sector $2\pi -\theta_0 < {\rm arg \; } z <2
\pi$ of the unphysical sheet. The necessary construction is given in section 2.
Then in sections 3 and 4 we show that our methods lead to certain resonance-free
regions which are adjacent to part of the real $z$-axis.
Finally, in section 5, we discuss the numerical computation of resonances
lying in the complement of our resonance-free regions. 
 \section{Analytic continuation}
\label{sec:2} The method which we develop in this section for continuing
$\psi(x,z)$ and $\psi^{'}(x,z)$ analytically into im $z<0$ is based on the
integral equation by means of which (\ref{eq:1.4})  is proved \cite[sections
1.3 and 1.9]{MSPE89}. Thus we begin by writing (\ref{eq:1.1}) (with $\lambda =
z^2 \neq 0$ and $y=\psi$) as a first-order system in a standard way   by
defining \beq W=\frac{1}{2}e^{-ixz}  \left ( \begin{array}{cc} 1 & -i/z \\ 1 &
i/z  \end{array} \right )  \left ( \begin{array}{c} \psi\\ \psi^{'} \end{array}
\right ). \label{eq:2.1}
\enq
Then 
\beq
W^{'} = \left \{ \left  (\begin{array}{cc}
0 & 0 \\
0 & -2iz \end{array} 
\right ) 
+ Q \left ( \begin{array}{cc}
-1 & -1 \\
1 & 1 
\end{array}
\right ) \right \} W
,
\label{eq:2.2}
\enq
where
\beq
Q=\frac{1}{2}iq/z,
\label{eq:2.3}
\enq
and the corresponding integral equation is
\beq
W(x,z)=e_1 + \int_x^\infty Q(t) K(t-x,z) W(t,z)dt,
\label{eq:2.4}
\enq
where
\beq
e_1 = \left ( \begin{array}{c}
1 \\ 0 \end{array} \right ),\;\;\;\;
K(s,z) = \left (
\begin{array}{cc}
1 & 1 \\
-e^{2isz} & -e^{2isz}
\end{array}
\right ).
\label{eq:2.5}
\enq
Iteration of (\ref{eq:2.4}) gives
\beq
W(x,z)=e_1+\sum_1^\infty W_n (x,z),
\label{eq:2.6}
\enq
where
\beq
W_n(x,z)=\int_x^\infty Q(t)K(t-x,z)W_{n-1}(t,z)dt
\label{eq:2.7}
\enq
and $W_0(x,z)=e_1$, provided of course that the infinite integral
converges. We note that, in terms of the components
of $W_n$ and $W_{n-1}$, (\ref{eq:2.7}) is
\beq
\left ( 
\begin{array}{c} u_n \\ v_n \end{array} \right ) (x,z)=
\int_x^\infty Q(t) \{ u_{n-1}(t,z) + v_{n-1}(t,z) \}
\left ( \begin{array}{c}
1 \\ - e^{2i(t-x)z} \end{array} \right ) dt.
\label{eq:2.8}
\enq
Also, the transformation (\ref{eq:2.1}) back to (\ref{eq:1.1}) gives
\beq
\psi=e^{ixz} \left ( 1 + \sum_1^\infty \{ u_n(x,z) + v_n(x,z) \} \right ),
\label{eq:2.9}
\enq
\beq
\psi^{'}=ize^{ixz} \left ( 1 + \sum_1^\infty \{ u_n(x,z) - v_n(x,z) \} \right ).
\label{eq:2.10}
\enq
In what follows, we write
\beq
\mid W_n \mid = \mid u_n \mid + \mid v_n \mid.
\label{eq:2.11}
\enq
\par
We now introduce the more detailed conditions on $q$ that we require.
\begin{condition}
We suppose that the real-valued function 
 $q(x)$ can be extended into a 
sector $S$ of the complex $\xi-$plane  as an analytic function
$q(\xi)$ as follows.
\begin{enumerate}
\item
$q(\xi)$ is regular in a sector $S$ defined by $-\theta_0 < {\rm arg \;} \xi
< \theta_0$ and $\xi \neq 0$, with some  $\theta_0$ such that $0 < \theta_0 <
\pi$.
\item
There are constants $\gamma$ $(>1)$ and $k$ such that
\beq
\mid q(\xi) \mid \leq  k \mid \xi \mid ^{- \gamma}
\label{eq:2.12}
\enq
as $\mid \xi \mid \rightarrow \infty$ and $\xi \in S$.
\end{enumerate}
\end{condition}
\par
Our method involves extending also the definition 
of $W(x,z)$ into the complex $\xi$-plane. To do this, we write $t=x+s$ in
(\ref{eq:2.7})  and consider the iterative definition
\beq
W_n(\xi,z)=\int_0^\infty Q(\xi+s)K(s,z)W_{n-1}(\xi+s,z)ds
\label{eq:2.13}
\enq
for $\xi \in S$, with $W_0(\xi,z)=e_1$.
In the following lemma we give a simple estimate for the size of $W_n$ 
in order to deal with the convergence of the infinite integral in
({\ref{eq:2.13}).
\begin{lemma}
Let $q$ satisfy Condition 2.1.
Let $\xi \in S$ and let im $z>0$ in (\ref{eq:2.5}) and (\ref{eq:2.13}).
Then for $n \geq 0$
\beq
\mid W_n(\xi,z) \mid \leq \frac{1}{n!}
\left ( 2 \int_0^\infty \mid Q(\xi + s )\mid ds \right )^n,
\label{eq:2.14}
\enq
and $W_n(\xi,z)$ is a regular function of $\xi$ in $S$.
\end{lemma}
{\bf Proof.}
We note that the infinite integral in (\ref{eq:2.14}) converges because of
(\ref{eq:2.12}). The lemma is clearly true when $n=0$ and, proceeding by
induction on $n$, we use the form of     (\ref{eq:2.8}) which
corresponds to (\ref{eq:2.13}).
By (\ref{eq:2.11}) and (\ref{eq:2.14}) (with $n-1$), this gives
\begin{eqnarray*}
\mid W_n(\xi,z) \mid & \leq &2^n \int_0^\infty \mid Q(\xi+s ) \mid
\frac{1}{(n-1)!}
\left (\int_0^\infty \mid Q(\xi+s + \sigma) \mid d \sigma \right )^{n-1} ds \\
&=&\frac{2^n }{(n-1)!}\int_0^\infty \mid Q(\xi+s ) \mid
\left (\int_s^\infty \mid Q(\xi + \sigma) \mid d \sigma \right )^{n-1} ds,
\end{eqnarray*}
from which (\ref{eq:2.14}) follows.
\par
To deal with the regularity of the $W_n$, we note that 
 (\ref{eq:2.12}) and (\ref{eq:2.14}) imply that the
infinite integral in (\ref{eq:2.13}) converges uniformly with respect to $\xi$
in any closed bounded region $S_1 \subset S$. Thus the regularity in $S$
of $W_n$ follows from that of $W_{n-1}$, and the lemma is proved.  
\par
The next step is to re-write (\ref{eq:2.13}) in a form which does not require
im $z>0$ and which therefore provides the analytic continuation of
$W_n(\xi,z)$ (as a function of $z$) into the lower half of
the $z-$plane. At this stage we restrict $\xi$ so that re $\xi\geq 0$,
the reason being given in the proof of the following theorem. Ultimately we
specialise $\xi$
to be the positive real variable $x$. 
 \begin{theorem}
Let $q$ satisfy Condition 2.1. Then, for all $\xi$
and $z$  in $S$ with re $\xi \geq 0$,
\beq
W_n(\xi,z) = \frac{1}{z} \int_0^\infty Q(\xi+\frac{s}{z})K(s,1)
W_{n-1}(\xi+\frac{s}{z},z) ds,
\label{eq:2.15}
\enq
and the series
\beq
W(\xi,z)=e_1+\sum_1^\infty W_n(\xi,z)
\label{eq:2.16}
\enq
defines a regular function of $z$ in $S$ which, when $\xi=x$, continues to
satisfy the differential equation (\ref{eq:2.2}).
\end{theorem}
{\bf Proof. }
We suppose first that im $z>0$, so that (\ref{eq:2.13}) holds. We consider the
contour integral
\beq
\int_C Q(\xi + \frac{\eta}{z})K(\eta,1)W_{n-1}(\xi + \frac{\eta}{z},z) d \eta
\label{eq:2.17}
\enq
where $C$ is the closed contour in the complex plane formed by the positive
real axis, the line through $z$ from $0$ to $\infty$,
and the smaller part of the circle  $\mid \eta \mid = R$.
The assumption that re $\xi \geq 0$ guarantees that the point $\xi+
\eta/z$ lies in $S$, and therefore the integrand in
(\ref{eq:2.17}) is defined as a regular function of $\eta$ within and on $C$.
Then, by Cauchy's Theorem, the value of (\ref{eq:2.17}) is zero.
Thus (\ref{eq:2.15}) follows from (\ref{eq:2.13}) when $R\rightarrow \infty$,
provided that the contribution to (\ref{eq:2.17}) from $\mid \eta \mid =R$
tends to zero.
\par
By (\ref{eq:2.3}), (\ref{eq:2.5}) and (\ref{eq:2.14}), this contribution does
not exceed in modulus
\beq
({\rm const.}) R \int_0^{ {\rm arg \; } z} \mid \xi + \frac{\eta}{z} \mid
^{-\gamma} \left ( \int_0^\infty \mid \xi + \frac{\eta}{z} + s \mid ^{-\gamma}
ds \right )^{n-1} d \theta
\label{eq:2.18}
\enq
in which $\eta = R e^{i \theta}$.
Since  $0 < {\rm arg \;} z < \theta_0 \; (<\pi)$, it is easy to check that
\begin{displaymath}
\mid X +i Y+ \frac{\eta}{z} \mid^2 \;\geq \frac{1}{2}(1-\mid \cos \theta_0 \mid
)(X^2 + R^2/\mid z \mid^2)-\kappa Y^2
\end{displaymath}
where $X= {\rm re \;} \xi + s,\;\;Y={\rm im }\; \xi$ and $\kappa = (1+3
\cos^2\theta_0)/\{(1-\mid \cos \theta_0 \mid )(1+3 \mid \cos \theta_0 \mid ) \}$.
Then, since $\gamma >1$, the $\theta-$integrand in (\ref{eq:2.18}) is
$O(R^{-\gamma -(n-1)(\gamma -1)})$, and hence (\ref{eq:2.18}) tends to zero as
$R \rightarrow\infty$ for all $n \geq1$. This proves (\ref{eq:2.15}) for im
$z>0$. \par We turn now to im $z \leq 0$ and we show that (\ref{eq:2.15})
continues to provide an iterative definition of the $W_n(\xi,z)$ as regular
functions of $z$.
We note that, when re $\xi \geq 0$ and im $z \leq 0$, the point $\xi + s/z$ in
(\ref{eq:2.15}) continues to lie in $S$.
An induction argument similar to that used for (\ref{eq:2.14}) shows that
\begin{equation}
\mid W_n (\xi,z) \mid \leq \frac{1}{n!} 
\left( \frac{2}{ \mid z \mid }\int_0^\infty \mid Q( \xi + \frac{s}{z}) \mid
ds \right)^n . \label{eq:2.19}
\end{equation}
Again, as for (\ref{eq:2.13}), the infinite integral in (\ref{eq:2.15})
converges uniformly with respect to $z$ in any closed bounded region $S_1
\subset S$, by (\ref{eq:2.12}).
Hence each $W_n(\xi,z)$ is a regular function of $z$ in $S$.
Further, (\ref{eq:2.19}) also guarantees the uniform convergence of the series
(\ref{eq:2.16}) with respect to $z$ in $S_1$,
and hence $W(\xi,z)$ is also a regular function of $z$ in $S$.
\par
Finally, we show that $W(\xi,z)$ satisfies (\ref{eq:2.2}) in the more general
form with $\xi$ in place of $x$. In (\ref{eq:2.15}), we sum for $n$ going from
$1$ to $\infty$ and we write $s=z(t-\xi)$ to obtain
\begin{equation}
W(\xi,z)=e_1 + \int_\xi^\infty Q(t) K( z (t-\xi),1)W(t,z) dt,
\label{eq:2.20}
\end{equation}
where $\infty$ denotes the point at infinity on the line through $\xi$ in the
direction of the vector $1/z$. The interchange of integration and summation
involved in (\ref{eq:2.20}) is justified by means of (\ref{eq:2.19}).
Differentiation of (\ref{eq:2.20}) with respect to $\xi$ now recovers
(\ref{eq:2.2}) with $\xi$ in place of $x$, and the proof of the theorem is
complete.
\par
We note that, in terms of the components of $W_n$ and $W_{n-1}$,
(\ref{eq:2.15}) is
\begin{equation}
\left( \begin{array} {c}
u_n \\v_n 
\end{array}
\right )
(\xi,z) = \frac{1}{2}i z^{-2} \int_0^\infty q(\xi+\frac{s}{z} )
(u_{n-1}+v_{n-1})( \xi + \frac{s}{z},z ) \left (
\begin{array}{c}
1 \\ -e^{2is} 
\end{array} 
\right) ds
\label{eq:2.21}
\end{equation}
corresponding to (\ref{eq:2.8}), and we have used (\ref{eq:2.3}).
Here (\ref{eq:2.21}) is valid for all $\xi$ and $z$ in $S$ with re $\xi \geq
0$, and then (\ref{eq:2.9}) and (\ref{eq:2.10}) provide the desired analytic
continuation of $\psi(x,z)$ and $\psi^{'}(x,z)$ into
the lower half of the sector $S$.
\section{Resonance-free regions}
The basic result on non-resonance which follows from (\ref{eq:2.9}),
(\ref{eq:2.10}) and (\ref{eq:2.21}) is given in the next theorem.
\begin{theorem}
Let $z \in S$ with {\rm im} $z<0$ and $z \neq i \cot \alpha$.
Let 
\begin{equation}
\int_0^\infty \mid q (s/z) \mid ds< \mid z \mid ^2 \log (1 + \delta ^{-1}),
\label{eq:3.1}
\end{equation}
where
\begin{equation}
\delta = \left \{
\begin{array}{cc}
1 & (0 \leq \alpha \leq \pi/2 \\
\mid \cos \alpha - i z \sin \alpha \mid / \mid \cos \alpha + i z \sin \alpha
\mid &(\pi/2 < \alpha < \pi ).
\end{array}
\right .
\label{eq:3.2}
\end{equation}
Then $\Psi(0,z) \neq 0$ and $z$ is not a resonance.
\end{theorem}
{\bf Proof.}
 By (\ref{eq:1.6}), (\ref{eq:2.9}) and (\ref{eq:2.10}), we have
\begin{equation}
\Psi(z)= ( \cos \alpha + i z \sin \alpha ) 
\left ( 1 + \sum_1^\infty \{ u_n(0,z) + Z v_n(0,z) \} \right ),
\label{eq:3.3a}
\end{equation}
where $Z= (\cos \alpha - i z \sin \alpha ) / (\cos \alpha + i z \sin \alpha)$.
It is easy to check that, since im $z<0$, $\mid Z \mid \leq 1$ for
$0 \leq \alpha \leq \pi/2$ and $\mid Z\mid > 1$
for $\pi/2 < \alpha < \pi$. Hence, with $\delta$ as in (\ref{eq:3.2}),
\begin{eqnarray*}
\mid \sum_1^\infty \{ u_n(0,z) + Zv_n(0,z) \} \mid &\leq& \delta
\sum_1^\infty \left (
\mid u_n(0,z) \mid + \mid v_n (0,z ) \mid \right ) \\
&\leq& \delta \left \{ \exp \left ( \mid z \mid ^{-2} \int_0^\infty \mid q
(s/z) \mid ds \right ) -1 \right \}
\end{eqnarray*}
by (\ref{eq:2.19}) and (\ref{eq:2.21}). It now follows from (\ref{eq:3.3a})
that $\Psi(z)$ is non-zero if $z \neq i \cot \alpha$ and
\begin{displaymath}
\exp \left ( \mid z \mid ^{-2} \int _0^\infty \mid q (s/z) \mid ds \right
){-1} < \delta^{-1},
\end{displaymath}
and the latter is guaranteed by (\ref{eq:3.1}).
\par
Let us note that, with the change of variable $s=\mid z \mid t$,
(\ref{eq:3.1}) can be written as
\begin{equation}
\int_0^\infty \mid q (t e^{i\theta})\mid dt < \mid z \mid \log (1+\delta^{-1}
\label{eq:3.3})
\end{equation}
where $\theta=2\pi-{\rm arg \;} z >0$. The condition (\ref{eq:3.3}) defines a
region of the complex plane within which there are no resonances,
and the nature of this resonance-free region depends on the nature of $q$.
Before we turn to detailed examples, we give one general property
of resonance-free regions which is a consequence of (\ref{eq:3.3}).
\begin{cor}
There are real numbers $R_1\;(>0)$ and $\theta_1$ $(0 < \theta_1 < \pi)$ such
that the sectorial region $\mid z \mid \geq R_1$,
$2 \pi -\theta_1  <  {\rm arg \;} z < 2\pi$ is resonance-free.
\end{cor}
{\bf Proof.} Suppose first that $0\leq \alpha \leq \pi/2$, so that
$\delta=1$ in (\ref{eq:3.2}). We choose $R_1$ so that
\begin{equation}
R_1 > (\log 2)^{-1} \int_0^\infty \mid q(t) \mid dt. \label{eq:3.5a}
\end{equation}
Then, by continuity in $\theta$, we have
\begin{equation}
\int_0^\infty \mid q(t e^{i\theta})\mid dt < R_1 \log 2
\label{eq:3.6a}
\end{equation}
for $\theta$ in some range $(0,\theta_1)$ with $\theta_1>0$.
Hence (\ref{eq:3.3}) holds for $\mid z \mid \geq R_1$, and the corollary is
proved for this range of $\alpha$.
\par
Next suppose that $\pi/2 < \alpha < \pi$.
Then, with $\delta$ as in (\ref{eq:3.2})
and $\cot \alpha < 0$, it is easy to check that
\begin{displaymath}
\delta \leq (1 + \sin \theta )/ \mid \cos \theta \mid
\end{displaymath}
where again  arg $z= 2 \pi -\theta$. We choose $R_1$ as in (\ref{eq:3.5a})
but, in place of (\ref{eq:3.6a}), we can say that
\begin{displaymath}
\int_0^\infty \mid q(te^{i \theta})\mid dt < R_1 \log \left (
1 + \frac{ \mid \cos \theta \mid }{1+ \sin \theta} \right ) < R_1 \log ( 1 +
\delta^{-1})
\end{displaymath}
for $\theta$ in some range $(0,\theta_1)$ with $\theta_1 >0$.
Hence (\ref{eq:3.3}) again holds for $\mid z \mid \geq R_1$, as required.
\par
Corollary 3.2 provides theoretical support for an observation by Aslanyan and
Davies \cite{AD98} concerning the numerical computation of resonances for the
potential
\begin{equation}
q(x) = x^2 \exp \left (-\epsilon x^2 \right ),
\label{eq:3.7a}
\end{equation}
where $\epsilon \;(>0)$ is a small parameter. In \cite[p. 16 and Table
10]{AD98} it is noted that there are resonances  very close to the positive
real axis but, at a certain point, they turn sharply away into the lower half
plane. Now (\ref{eq:3.7a}) satisfies Condition 2.1 with $\theta_0 < \pi/4$
(see also Example 4.6 below), and the existence of the sectorial region in
Corollary 3.2
precludes as a general  feature the occurrence of resonances close to the
positive real axis beyond a certain distance from the origin.
\par
Corollary 3.2 can also be related to \cite[Theorem 1]{MSPE98a} concerning the
localization of spectral concentration points to a bounded interval on the
real spectral axis (see also \cite[section 2]{BE99}). Insofar as spectral
concentration is associated with resonances located near to the real axis,
Corollary 3.2 provides another proof that spectral concentration points are
confined to a bounded interval for a class of potentials satisfying
(\ref{eq:1.3}).
\par
In the Dirichlet case $\alpha=0$ of (\ref{eq:1.2}), there are additional non-resonance results like Theorem 3.1 and Corollary 3.2 but with, in the corollary, the vertex of the sector at the origin.
\begin{theorem}
Let $\alpha=0$ in (\ref{eq:1.2}) and let $q$ satisfy Condition 2.1 with
\beq
\gamma>2 \label{eq:3.8a}
\enq
in (\ref{eq:2.12}). Let $z \in S$ with {\rm im }$z<0$, and let
\beq
\int_0^\infty s \mid q(s/z) \mid ds < \mid z\mid^2 \log 2.
\label{eq:3.9a}
\enq
Then $\psi(0,z)\neq 0$ and $z$ is not a resonance.
\end{theorem}
{\bf Proof }.  We note that (\ref{eq:3.8a}) guarantees the convergence of the integral in (\ref{eq:3.9a}). In (\ref{eq:2.21}), we use the inequality
$\mid 1 - e^{2is} \mid \leq 2s$ to obtain
\begin{displaymath}
\mid (u_n + v_n ) (\xi,z) \mid \leq \mid z \mid ^{-2} \int_0^\infty s \mid q( \xi + \frac{s}{z}) \mid \mid (u_{n-1} + v_{n-1})(\xi + \frac{s}{z},z) \mid ds.
\end{displaymath}
Then, as for (\ref{eq:2.19}), an induction argument gives
\begin{displaymath}
\mid (u_n + v_n ) (\xi,z) \mid \leq \frac{1}{n!}
\left ( \frac{1}{\mid z \mid^2} \mid \int_0^\infty s \mid q( \xi + \frac{s}{z}) \mid 
ds \right )^n.
\end{displaymath}
This inequality is used in (\ref{eq:2.9}) and (\ref{eq:1.6}) (with $\alpha =0$),
and the theorem follows from (\ref{eq:3.9a}) in the same way as Theorem 3.1 followed from (\ref{eq:3.1}).
\par
As for (\ref{eq:3.3}), the change of variable $s=\mid z\mid t $ in (\ref{eq:3.9a}) leads to
\beq
\int_0^\infty t \mid q (t e^{i\theta})\mid dt < \log 2
\label{eq:3.10}
\enq
and this in turn leads immediately to the next corollary.
\begin{cor}
Let $\alpha=0$ in (\ref{eq:1.2}) and, in addition to (\ref{eq:3.8a}), let
\beq
\int_0^\infty t\mid q(t) \mid dt < \log 2.
\label{eq:3.11}
\enq
Then there is a real number $\theta_1\;(0<\theta_1 < \pi)$ such that the sector $\mid z \mid >0,$ $2\pi - \theta_1\; < {\rm arg } \;z < 2\pi$ is resonance-free.
\end{cor}
\par
The condition (\ref{eq:3.11}) can be  related to the condition
\beq
\int_0^\infty t \mid q(t) \mid dt < 0.1735
\label{eq:3.12}
\enq
\cite[(2.16)]{MSPE98b} which is shown in \cite{MSPE98b} (by quite different methods) to imply the absence of any spectral concentration points on the positive spectral axis $(0,\infty)$.
The smaller the value of the integral in (\ref{eq:3.11}), the larger $\theta_1$ can be, and the further away from the real axis are any resonances pushed. Thus, in the case of (\ref{eq:3.12}), any resonances are too far from the real axis to produce spectral concentration \cite[section 3(iv)]{MSPE98b}.

 \section{Examples}
We consider now some examples of
$q$ which show in more detail the type of region that arises from
(\ref{eq:3.3}). We keep to the case $0 \leq \alpha \leq \pi/2$ for which
$\delta=1$ in (\ref{eq:3.2}): in the other case, $\delta \rightarrow 1$ as
$\mid z\mid \rightarrow \infty$ and the regions are asymptotically similar for
large $\mid z \mid$.
\subsection{Example $q(x)=c(x+a)^{-\gamma}$} where $\gamma>1$, $c$ and $a$ are real and
$a>0$. In Condition 2.1, we take $q(\xi)=c(\xi+a)^{-\gamma}$ with, if $\gamma$
is not an integer, a cut in the $\xi-$plane from $-a$ to $-\infty$ along the
real axis. Thus we can take $\theta_0=\pi$. The integral in  (\ref{eq:3.3}) is
now 
\begin{equation} \mid c \mid\int_0^\infty \mid t + a e^{-i \theta} \mid^{-\gamma} dt=
\mid c \mid\int_0^\infty(t^2+2at \cos \theta + a^2 )^{-\gamma/2}dt = I(\theta)
\label{eq:3.4}
 \end{equation}
say. Hence $I(\theta)$ increases from $\mid c \mid a^{-\gamma+1}/(\gamma-1)$ to $\infty$
as $\theta$ increases from $0$ to $\pi$, and (\ref{eq:3.3}) becomes
\begin{equation}
\mid z \;\mid > I(\theta)/\log 2.
\label{eq:3.5}
\end{equation}
Thus we have a resonance-free region which lies in the lower half of the the
complex plane,
bounded by a    curve which starts at the point $\mid c \mid a^{-\gamma+1}/(\gamma-1)$
on the real axis and recedes from the origin as $\theta\;(=-{\rm arg \; } z)$
increases from $0$ to $\pi$.
The region is of course on the side of the curve remote from the origin. When
$\gamma=2$ in particular, the integration in (\ref{eq:3.4}) can be performed
and (\ref{eq:3.5}) becomes
\begin{displaymath}
\mid z \mid \;> \frac{\mid c \mid}{a \log 2} \frac{\theta}{\sin \theta}\;\;\;(\theta =
-{\rm arg \; } z).
\end{displaymath}
Thus the boundary curve in this case is asymptotic from above to the line im $z=-\mid c \mid \pi / (a \log 2)$ as $\theta \rightarrow \pi.$
\subsection{Example $q(x)=c(x^n + a^n)^{-\gamma}$}
where $n\;(\geq 2)$ is an integer, $n\gamma >1$, $c$ and  $a$ are real and $a>0$.
This is similar to Example 3.1 but now $\theta_0=\pi/n$.
The integrand\ in (\ref{eq:3.4}) is replaced by $(t^{2n} + 2 a^nt^n\cos n 
\theta + a^{2n})^{-\gamma/2}$, and $I(\theta)$ increases to $\infty$
as $\theta \rightarrow \pi/n$.
\par
In the case when $n=2$ and $\gamma=2$, we find that
\begin{displaymath}
I(\theta)= ( \pi \mid c \mid / 4 a^3 )\sec \theta,
\end{displaymath}
and the resonance-free region (\ref{eq:3.5}) is the quadrant
\begin{displaymath}
{\rm re}\; z > \pi \mid c \mid / (4 a^3 \log 2 ),\;\;\;\; {\rm im }\; z <0.
\end{displaymath}
Also in this case, the Dirichlet condition (\ref{eq:3.10}) gives
\begin{displaymath}
2\theta/ \sin 2 \theta < 2 ( a^2/ \mid c \mid )\log 2.
\end{displaymath}
Thus, on the assumption that $\mid c \mid < 2a^2 \log 2 $, this being (\ref{eq:3.11}), the value of $\theta_1$ in Corollary 3.4 is the solution of
\begin{displaymath}
2 \theta_1 / \sin 2 \theta_1 = 2 ( a^2/\mid c \mid )\log 2.
\end{displaymath}

\subsection{Example $q(x)=c\{(x-w)(x-\bar{w})\}^{-\gamma}$}
where $2 \gamma >\;1$, $w \neq 0$ and $0 < {\rm arg \; } w < \pi$.
This again is similar. Here $\theta_0={\rm arg \; } w \;(=\phi$, say ),
and the integrand in (\ref{eq:3.4}) is replaced by
\begin{displaymath}
\{ (t^2-\mid w \mid^2 )^2-4 \mid w \mid t ( t - \mid w \mid )^2
\cos \phi \cos \theta + 4 \mid w \mid ^2 t^2 ( \cos \theta -\cos \phi )^2
\}^{-\gamma/2}.
\end{displaymath}
Again $I(\theta) \rightarrow \infty$ as $\theta \rightarrow \phi$ because
we approach a singularity at $t=\mid w \mid$ in $\{...\}^{-\gamma/2}$.
However, $I(\theta)$ is not necessarily monotonic unless $\cos \phi \leq 0$.
\par
We conclude this group of examples by noting that similar remarks apply when
$q$ is a product of terms already considered with differing values of $a$,
$w$, and $\gamma$ and, indeed, when $q$ is a ratio of two such products.
We give one example of this more general type for future reference in Section 5.
\subsection{Example $q(x)=c(x-1)/(x+1)^4.$}
Here $I(\theta)$ in (\ref{eq:3.4}) and (\ref{eq:3.5}) is replaced  by
\beq
I(\theta) = \mid c\mid \int_0^\infty \mid t e^{i \theta} -1 \mid /(t^2+2t \cos \theta +1)^2 dt. \label{eq:4.3a}
\enq
Now the boundary curve of the resonance-free region(\ref{eq:3.5}) starts at the point $0.36 \mid c \mid$ on the real axis and, since $I(\theta) \sim ({\rm const.}) (\sin \theta )^{-3}$ when $\theta \rightarrow \pi$, the curve is asymptotically like
$\mid \rm{ re}\;  z  \mid^2=({\rm const.})\mid {\rm im \;} z \mid^3$
(see also Figure 1 below). \par
Next, we turn to examples with exponential decay which are also covered by
Condition 2.1.
\subsection{Example $q(x)=2e^{-ax}\sin x \;\;(a>0).$}
In Condition 2.1, we take
\begin{equation}
q(\xi)=i \left (e^{-(a+i)\xi}- e^{-(a-i)\xi}\right ).
\label{eq:3.6}
\end{equation}
Since $\mid e^{-(a\pm i)\xi} \mid = e^{-a {\rm re \; } \xi \pm {\rm im \;}\xi}$,
(\ref{eq:2.12}) is certainly satisfied if
\begin{equation}
\theta_0 < \tan^{-1}a.
\label{eq:3.7}
\end{equation}
By (\ref{eq:3.6}), the left-hand side of (\ref{eq:3.3}) does not exceed
\begin{displaymath}
2 \int_0^\infty \exp \left \{ \left (
-a \cos \theta + \sin \theta \right ) t \right \} dt=
2(a \cos \theta - \sin \theta)^{-1}.
\end{displaymath}
Hence (\ref{eq:3.3}) holds if $\mid z \mid (a \cos \theta - \sin \theta )>
2/\log 2$, or
\begin{displaymath}
a ({\rm re \;} z) +({\rm im \; }z)>  2/\log 2.
\end{displaymath}
Since $\theta_0$ can be arbitrarily near to $\tan^{-1} a $ in (\ref{eq:3.7}),
we therefore have a resonance-free region in the lower half of the complex
plane lying to the right of the line through the point $2/(a\log 2)$ on the 
real axis and with gradient $-a$.
\par
We observe  that independent support for this gradient $-a$ is provided by the
quite different analytic continuation
method developed in \cite[Prop. 2.1]{BEM99}.
This latter method constructs   the analytic continuation of $\Psi(z)$
into the whole of im $z<0$ except for poles at the points
\begin{displaymath}
 -\frac{1}{2} \left ( \nu + m a i \right )\;\;\;\;\left ( \mid \nu \mid\leq
m,\;\;m=1,2,... \right ),
\end{displaymath}
$\nu$ being an integer. Thus the line through the origin with the same
gradient $-a$ delineates a pole-free region for
$\Psi(z)$ within which $\Psi(z)$ is regular.
The methods in \cite{BEM99} do not however lead readily to resonance-free
regions.
\subsection{Example $q(x)=cx^m\exp(-x^n)$ \tiny \normalsize}
where $m$ and $n$ are positive integers.
In Condition 2.1, we take
\begin{displaymath}
q(\xi)=c\xi^m \exp (-\xi^n),
\end{displaymath}
and (\ref{eq:2.12}) is certainly satisfied if
\begin{equation}
\theta_0 < \pi/2n.
\label{eq:3.8}
\end{equation}
Now the left-hand side of (\ref{eq:3.3}) is
\begin{displaymath}
\mid c\mid\int_0^\infty t^m \exp \left ( -t^n \cos n \theta \right ) dt = \mid c\mid\left (\cos n
\theta \right )^{-(m+1)/n} I,
\end{displaymath}
where $I=\int_0^\infty u^m \exp (-u^n) du$. Hence (\ref{eq:3.3}) holds if
\begin{equation}
\mid z \mid ( \cos n \theta )^{(m+1)/n} >\mid c \mid I/\log 2.
\label{eq:3.9}
\end{equation}
In (\ref{eq:3.8}), $\theta_0$ can be arbitrarily near to $\pi/2n$ and hence,
in (\ref{eq:3.9}), we can let $\theta$ increase from $0$ to $\pi/2n$.
Thus (\ref{eq:3.9}) defines a region in the lower half plane whose boundary
starts at the point $\mid c \mid I/\log 2$ on the real axis and recedes to infinity as
$\theta \rightarrow \pi/2n$.
\par  
A typical example of (\ref{eq:3.9}) is when $m=0$ and $n=2$, in which case
the boundary is the part of the rectangular hyperbola $X^2-Y^2=(\mid c \mid
I/\log 2)^2$ which lies in the fourth quadrant of the $(X,Y)-$plane and
$z=X+iY$. Again the resonance-free region lies on the side of the hyperbola
remote from the origin.
Independent support for the nature of this boundary is provided by the findings
of Siedentop \cite{HS87a} and Froese \cite{F97}. In \cite{HS87a}
(where $c=-1$), the first
few resonances found computationally are already near to, but below, the line
arg $z=-\pi/4$ (see also \cite[Example 6.4]{BEM99}) while, in \cite{F97},
the resonances are shown to be asymptotically near to this same line.
\section{Computational resonance-finding}
We turn now to the numerical computation of resonances for explicit $q$, such as those in section 4, which satisfy Condition 2.1. One possible direct method
is to compute the $u_n$ and $v_n$ $(1 \leq n \leq N)$ recursively in (\ref{eq:2.21}) and substitute the results into (\ref{eq:2.9}) and (\ref{eq:2.10}), the infinite series being truncated at $N$ with an  error term.
Then a zero-finding algorithm would be applied to the resulting approximation to $\Psi$ in (\ref{eq:1.6}).
A similar procedure was applied successfully to the formulae  for $u_n$ and $v_n$ in \cite{BEM99} when $q$ has exponential decay.
However, in our present situation of power decay, it has proved difficult to use (\ref{eq:2.21}) when $N \geq 2$, repeated integration being
involved, and in addition the error term for $N=1$ is not small.
\par
Instead, we have computed resonances by the method of complex scaling. We refer to Simon \cite{BS78} for a discussion of this method in relation to resonances and to Agmon \cite{Ag98} for a recent definitive account in a very general setting.
The method of complex scaling is closely associated with (\ref{eq:2.20}) and (\ref{eq:2.9}) and, in fact, our approach in section 2 provides an independent justification of the validity of this method for (\ref{eq:1.1}), as we now describe.
\par
The transformation of (\ref{eq:2.20}) back to $\psi(\xi, z)$ via (\ref{eq:2.1})
(with $\xi$ in place of $x$) gives
\begin{displaymath}
d^2\psi/d \xi^2 + \{ z^2-q(\xi)\} \psi =0
\end{displaymath}
corresponding to (\ref{eq:1.1}).
With $\xi$ in polar form
$\xi=r \exp(i\phi)$ $(0<\phi <\theta_0)$, we therefore have
\beq
e^{-2i\phi} d^2\psi/dr^2 + \{ z^2 -q(r e^{i\phi})\}\psi=0
\label{eq:5.1}
\enq and, by (\ref{eq:2.9}) and (\ref{eq:2.19}),
\begin{displaymath}
\mid \psi (re^{i \phi},z) \mid =
[\exp \{ -r \mid z \mid \sin (\phi + {\rm arg}\; z ) \}]\{1 + o(1) \},
\end{displaymath}
where $o(1)$ refers to $r \rightarrow \infty$. It follows that $\psi(r e ^{i \phi},z)$ is an $L^2(0,\infty)$ solution of (\ref{eq:5.1}) if
\beq
2\pi -\phi < {\rm arg }\; z < 2 \pi.
\label{eq:5.2}
\enq
Thus the zeros of $\Psi(z)$ in (\ref{eq:1.6}) provide the eigenvalues $z^2$ of (\ref{eq:5.1}) on $0\leq r < \infty$ with the boundary condition
\begin{displaymath}
\psi(0,z)\cos \alpha + \psi^{'}(0,z) \sin \alpha =0
\end{displaymath}
at $r=0$. Here (\ref{eq:5.1}) is said to be obtained from (\ref{eq:1.1}) by complex scaling, the scaling factor being $e^{i \phi}$ \cite[section 5]{AD98}
\cite[section 3]{BS78}.

We have therefore applied a computational eigenvalue finder \cite{GM99} to (\ref{eq:5.1}) with a suitable value of $\phi$.
This locates eigenvalues $z^2$ and hence resonances in the sector (\ref{eq:5.2}). We focus the discussion of our computational findings now on Examples 4.1-4.4 since it is potentials with only power decay which are the main object of this paper.
\par
We   consider first
\beq
q(x)=c (x^2+1)^{-2}, \label{eq:5.3}
\enq
being the case $n=\gamma=2,\;a=1$ of Example 4.2.
Here $\theta_0=\pi/2$ but, if $\phi$ in (\ref{eq:5.1}) is close to $\pi/2$,
the code in \cite{GM99} reports unreliable results due to the sharp (but non-singular)
maximum of the scaled $\mid q(r e^{i \phi})\mid$ near to $r=1$.
Accordingly we have chosen $\phi=1.5$. We have found no resonances satisfying
(\ref{eq:5.2}) within the disk $\mid z \mid < 10$ when $c$ has the range of
values $-1,-5,-10,-15,-20$. This is certainly consistent with the
resonance-free quadrant \mbox{re $z > \pi \mid c \mid / ( 4 \log 2),$ im $z<0$}
in Example 4.2, but there remains the open question whether resonances occur
elsewhere in im  $z<0$. \par A similar example, but with  a  higher
singularity located nearer to the real axis in the complex plane, is \beq
q(x)=c(x^6+1)^{-20}. \label{eq:5.4} \enq Despite this extra feature, this
example also produces no spectral concentration and no resonances. Here
$\theta_0=\pi/6$ and we have chosen $\phi=0.5$.
The values of $c$ investigated were $-1,-5,-15,-20,-25,-30$.
The reason for choosing $c$ negative in (\ref{eq:5.3}) and (\ref{eq:5.4})
is to give $q(x)$ a negative minimum, a property which in exponentially 
decaying examples is often associated with spectral concentration and resonances
\cite[Section 6]{BEM99}.
\par
Next we consider Example 4.4
\begin{displaymath}
q(x)=c(x-1)/(x+1)^4,
\end{displaymath}
this time with $c>0$ to give the negative minimum (at $x=0$).
Here of course $\theta_0=\pi$ and we have chosen $\phi=3.0$.
There is one real point of spectral concentration when $c=35$ located at
$\lambda=0.26$, for which $\sqrt \lambda =0.51$, and we have tracked the
corresponding resonance for a range of values down to $c=0.5$. The resonance
broadly recedes from the real $z-$axis as $c$ decreases, and we give a
selection of these findings in Table 1. \begin{table}[htbp] \begin{center}
\begin{tabular}{||c|c|c ||} \hline  $c$ &  re$z$  &im  $z$  \\ \hline
  $35$ & 0.50  & -0.06    \\ \hline
 $25$ & 0.65  & -0.23    \\ \hline
 $15$ & 0.57  & -0.44    \\ \hline
 $10$ & 0.42  & -0.55    \\ \hline
 $4$ & 0.06  & -0.68    \\ \hline
 $3$ & -0.05  & -0.68    \\ \hline
 $2$ & -0.21  & -0.66    \\ \hline
 $1$ & -0.45  & -0.58    \\ \hline
 $0.5$ & -0.67  & -0.45    \\ \hline

 \end{tabular}
\end{center}
\caption{Resonance $z$}
\label{tab:3}
\end{table}
For small values of $c$, re $z$ appears to increase rapidly in the negative direction, but arg $z$ becomes too close to $\pi$ for the code \cite{GM99} to produce reliable values.
\par
 In order to gain an idea of how Table 1 relates to the resonance-free region given by (\ref{eq:3.5}) 
and (\ref{eq:4.3a}), we note that $I(\theta)$ contains a factor $\mid c
\mid$. Accordingly, we have applied a scaling factor $\mid c \mid ^{-1}$ to
both $I(\theta)$ and the values in Table 1. The result is Figure  1, in which the diamonds denote the scaled resonances from Table 1, and the dotted curve denotes the boundary curve scaled to $c=1$.  Figure 1 confirms the general nature of our theoretical result (\ref{eq:3.3}).
\par
We also mention that there are two additional similar strings of resonances: when $c=10$ for example, there are resonances at $z=-1.27-1.39i$ and $z=-5.05-5.31i$ in addition to the value in Table 1.
These additional resonances, however, lie  further from the resonance-free region than the resonance-string shown in Figure 1.
 \begin{figure}[htbp]
\centerline{
\epsfxsize=10cm
\epsffile {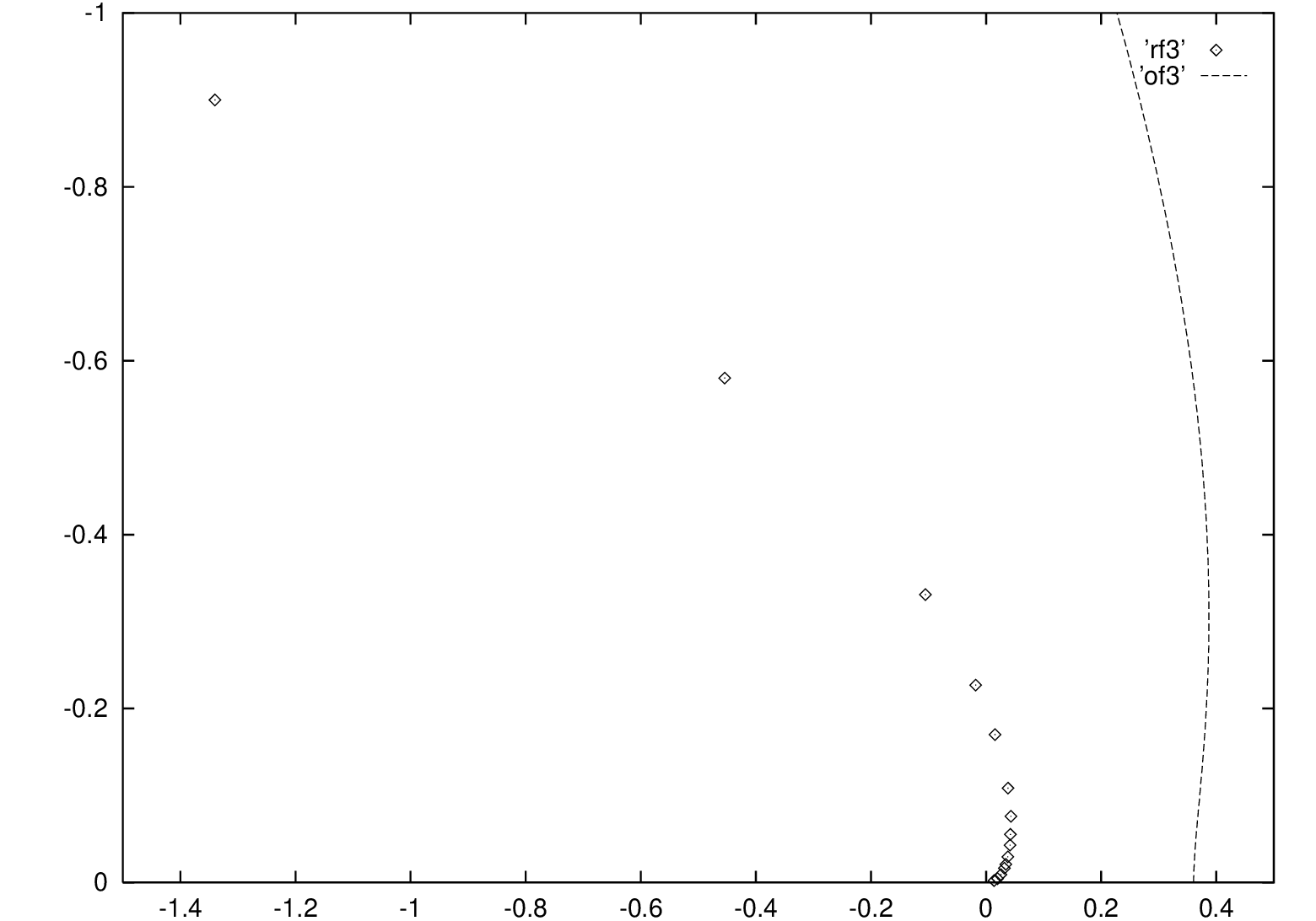}
}
\caption{}
\end{figure}

\par
Finally, we have also considered the example \begin{displaymath}
q(x)=c(x-1)/(x^4+1) \end{displaymath}
for which $\theta_0=\pi/4$ and we have taken $\phi=0.75$.
There is one real point of spectral concentration when $c=11$ located at $\lambda=0.15$, for which $\sqrt \lambda=0.39$.
We have tracked the corresponding resonance from $z=0.39 -0.03 i$ when $c=11$ as far as $z=0.73-0.59i$ when $c=3.2$.
For smaller $c$, the code \cite{GM99} again flags unreliability, but there is
a corresponding picture to Figure 1 to similarly confirm the theoretical result (\ref{eq:3.3}).

 \bibliographystyle{plain}
\bibliography{../bib_dir/bib1,../bib_dir/bibliography,../bib_dir/specon,../bib_dir/specon2,../bib_dir/Zettl,../bib_dir/help,../bib_dir/specon1,../bib_dir/analytic}  
\end{document}